\documentclass[oneside,english]{amsart}
\usepackage[T1]{fontenc}
\usepackage[latin9]{inputenc}
\usepackage{amsthm}
\usepackage{amssymb}
\usepackage{esint}

\makeatletter
\numberwithin{equation}{section}
\numberwithin{figure}{section}
\theoremstyle{plain}
\newtheorem{thm}{\protect\theoremname}
  \theoremstyle{plain}
  \newtheorem{lem}[thm]{\protect\lemmaname}

\makeatother

\usepackage{babel}
  \providecommand{\lemmaname}{Lemma}
\providecommand{\theoremname}{Theorem}

\begin{document}

\title{On Zeros of Fourier Transforms}

\author{Ruiming Zhang}
\begin{abstract}
In this work we verify the sufficiency of a Jensen's necessary and
sufficient condition for a class of genus 0 or 1 entire functions
to have only real zeros. They are Fourier transforms of even, positive,
indefinitely differentiable, and very fast decreasing functions. We
also apply our result to several important special functions in mathematics,
such as modified Bessel function $K_{iz}(a),\ a>0$ as a function
of variable $z$, Riemann Xi function $\Xi(z)$, and character Xi
function $\Xi(z;\chi)$ when $\chi$ is a real primitive non-principal
character satisfying $\varphi(u;\chi)\ge0$ on the real line, we prove
these entire functions have only real zeros.
\end{abstract}


\curraddr{College of Science\\
 Northwest A\&F University\\
 Yangling, Shaanxi 712100\\
 P. R. China.}


\email{ruimingzhang@yahoo.com}

\thanks{This work is partially supported by the National Natural Science
Foundation of China, grant No. 11371294. The author specially thanks
his colleague Dr. Bing He for helping him to check some of the computations.}

\maketitle

\section{Introduction}

One of fundamental questions in the theory of special functions is
whether an entire function has only real zeros. Some interesting examples
are even entire functions of genus $0$ or $1$, for example, Jackson's
$q$-Bessel function $J^{(2)}(z;q)$ and Bessel function of first
kind $J_{\nu}(z)$, they have only real zeros, and they are essentially
even entire functions of genus 0 and 1 respectively, \cite{Andrews,Ismail}.
The genus of an entire function $f(z)$ is an concept mostly related
to its Hadamard infinite product expansion and its order of growth.
For each $r>0$, let \cite{Levin}
\begin{equation}
\left\Vert f\right\Vert _{\infty,r}=\sup_{\left|z\right|\le r}\left|f(z)\right|=\max_{\left|z\right|=r}\left|f(z)\right|.\label{eq:1.1}
\end{equation}
Then the order $\rho$ of $f(z)$ can be defined by

\begin{equation}
\rho=\limsup_{r\rightarrow\infty}\frac{\log\left(\log\left\Vert f\right\Vert _{\infty,r}\right)}{\log(r)}.\label{eq:1.2}
\end{equation}
Given an arbitrary point $z_{0}\in\mathbb{C}$, the order $\rho$
of $f(z)$ can be computed from the values of its derivatives at $z_{0}$,
\begin{equation}
\rho=\left(1-\limsup_{n\to\infty}\frac{\log\left|f^{(n)}(z_{0})\right|}{n\log n}\right)^{-1}.\label{eq:1.3}
\end{equation}
Clearly, formula (\ref{eq:1.3}) shows that $f^{'}(z)$ has the same
order as $f(z)$. If the order $\rho$ of $f(z)$ is finite, then
it has the following Hadamard factorization

\begin{equation}
f(z)=z^{m}e^{p(z)}\prod_{n=1}^{\infty}\left(1-\frac{z}{z_{n}}\right)\exp\left(\frac{z}{z_{n}}+\ldots+\frac{1}{k}\left(\frac{z}{z_{n}}\right)^{k}\right),\label{eq:1.4}
\end{equation}
where $m$ is an nonnegative integer, $\left\{ z_{n}\right\} _{n=1}^{\infty}$
is the set of all nonzero roots of $f(z)$, $p(z)$ a polynomial of
degree $j$, and $k$ is the smallest nonnegative integer such that 

\begin{equation}
\sum_{n=1}^{\infty}\frac{1}{|z_{n}|^{k+1}}<\infty.\label{eq:1.5}
\end{equation}
The genus of $f(z)$ is defined as $g=\max\left\{ j,k\right\} $.
If the order of $f(z)$ is not an integer, then $g=\left\lfloor \rho\right\rfloor $
is the integer part of $\rho$, otherwise, $g$ may be $\rho-1$ or
$\rho$. 

In 1913 Jensen proved a set of necessary and sufficient conditions
for a class of genus 0 or 1 entire functions to have only real zeros,
\cite{Gasper1,Gasper2,Polya}. Gasper applied Jensen's conditions
to prove many important special functions that have only real zeros
by directly verifying positivities of certain sums and integrals,
\cite{Gasper1,Gasper2}. In this work we verify a Jensen inequality
for a class of functions under much stronger conditions. More specifically,
we are going to show the following: 
\begin{thm}
\label{thm:gasper} Let $\varphi(u)$ be nonnegative, even, and indefinitely
differentiable function that is not identically zero. If additionally
for each $n\in\mathbb{Z}^{+}$ there exists some $d,\,\delta>0$ such
that 
\begin{equation}
\varphi^{(n)}(u)=\mathcal{O}\left(\exp\left(-de^{\delta\left|u\right|}\right)\right)\label{eq:gasper}
\end{equation}
as $u\to\pm\infty$. Then the Fourier transform $\Phi(z)$ of $\varphi(u)$
has only real zeros.
\end{thm}
We just mention three applications of our results to special functions
here. For $\Re(z)>0$, the Bessel functions $K_{\nu}(z)$ has integral
representation \cite{Andrews,Erdelyi,Gasper2,Watson}
\begin{align}
K_{\nu}(z) & =\frac{1}{2}\int_{-\infty}^{\infty}\exp\left(-z\cosh u-\nu u\right)du.\label{eq:1.6}
\end{align}
Then for $a>0$ 
\begin{equation}
K_{iz}(a)=\frac{1}{2}\int_{-\infty}^{\infty}e^{-a\cosh u}e^{izu}du\label{eq:1.7}
\end{equation}
defines an entire function of variable $z$. Clearly, $\varphi(u)=e^{-a\cosh u}$
satisfies all the conditions of Theorem \ref{thm:gasper}. Thus the
entire function $K_{iz}(a)$ has only real zeros.

The Riemann Xi function is defined by \cite{Biane,Csordas,Edwards,Gasper1,Gasper2,Titchmarsh}
\begin{equation}
\Xi\left(s\right)=-\frac{1+4s^{2}}{8\pi^{\frac{1+2is}{4}}}\Gamma\left(\frac{1+2is}{4}\right)\zeta\left(\frac{1+2is}{2}\right),\label{eq:1.8}
\end{equation}
where $\Gamma(z)$ and $\zeta(z)$ are analytic continuations of Euler
gamma and Riemann zeta functions respectively. Then the entire function
$\Xi\left(s\right)$ satisfies \cite{Biane,Csordas,Edwards,Gasper1,Gasper2,Polya,Titchmarsh}
\begin{equation}
\Xi\left(s\right)=\int_{-\infty}^{\infty}\varphi(u)e^{isu}du,\quad s\in\mathbb{C},\label{eq:1.9}
\end{equation}
where 
\begin{align}
\varphi(u) & =2\pi\sum_{n=1}^{\infty}\left\{ 2\pi n^{4}e^{9u/2}-3n^{2}e^{5u/2}\right\} \exp\left(-n^{2}\pi e^{2u}\right),\ u\in\mathbb{R}.\label{eq:1.10}
\end{align}
It is known that $\varphi(u)$ is even, positive, indefinitely differentiable.
Furthermore, for any $\epsilon>0$ and integers $n\ge0$, it satisfies
\begin{align}
\varphi^{(n)}(u)= & \mathcal{O}\left(\exp\left\{ -(\pi-\epsilon)e^{2\left|u\right|}\right\} \right)\label{eq:1.11}
\end{align}
as $u\to\pm\infty$. Hence, it also satisfies the requirements of
our result Theorem \ref{thm:gasper}. Thus $\Xi\left(s\right)$ has
only real zeros, which means the Riemann hypothesis is valid.

Given a positive integer $q\ge2$, let $\chi(n)$ be a primitive real
character with respect to modulus $q$ with parity $a$, \cite{Davenport,Titchmarsh}
\begin{equation}
a=\begin{cases}
0,\  & \chi(-1)=1\\
1,\  & \chi(-1)=-1
\end{cases}.\label{eq:1.12}
\end{equation}
Define the character Xi function $\Xi\left(s;\chi\right)$ for $\chi$
of parity $a$ by \cite{Davenport} 
\begin{equation}
\Xi\left(s;\chi\right)=\frac{\Gamma\left(\frac{2a+1+2is}{4}\right)L\left(\frac{1}{2}+is,\chi\right)}{\left(\frac{\pi}{q}\right)^{\left(2a+1+2is\right)/4}},\label{eq:1.13}
\end{equation}
where $\Gamma(s)$ and $L\left(s,\chi\right)$ are the analytic continuations
of Euler gamma function and Dirichlet $L$-series for $\chi$ respectively.
Then the entire function $\Xi\left(s;\chi\right)$ satisfies 
\begin{equation}
\Xi\left(s;\chi\right)=\int_{-\infty}^{\infty}e^{isu}\varphi\left(u;\chi\right)du,\label{eq:1.14}
\end{equation}
where
\begin{equation}
\varphi\left(u;\chi\right)=2e^{w_{a}u}\sum_{n=1}^{\infty}n^{a}\chi(n)e^{-n^{2}\pi e^{2u}/q},\quad w_{a}=\frac{1+2a}{2},\label{eq:1.15}
\end{equation}
and $a$ is the parity of $\chi$. By the transformation formulas
\cite{Davenport} 
\begin{equation}
\sum_{n=-\infty}^{\infty}n^{a}\chi(n)e^{-n^{2}\pi u/m}=u^{-w_{a}}\sum_{n=-\infty}^{\infty}n^{a}\chi(n)e^{-n^{2}\pi/(mu)},\quad u>0\label{eq:1.16}
\end{equation}
we have
\begin{equation}
\varphi\left(-u;\chi\right)=\varphi\left(u;\chi\right),\quad u\in\mathbb{R}.\label{eq:1.17}
\end{equation}
It is clear that for each integer $n\ge0$ we have
\begin{equation}
\varphi^{(n)}\left(u;\chi\right)=\mathcal{O}\left(\exp\left(-de^{2|u|}\right)\right)\label{eq:1.18}
\end{equation}
as $u\to\pm\infty$ for any $0<d<\frac{\pi}{q}$. Hence for those
$\chi$ such that $\varphi\left(u;\chi\right)>0,\ u\in\mathbb{R}$,
$\varphi\left(u;\chi\right)$ satisfies all the conditions of Theorem
\ref{thm:gasper}. Thus the entire function $\Xi\left(s;\chi\right)$
has only real zeros, hence the generalized Riemann hypothesis holds
for this kind of Dirichlet $L$-series.

\section{Main Results}

The following theorem only states a subset of Jensen's results related
to the current work \cite{Gasper1,Polya} :
\begin{thm}
\label{thm:jensen}Let $\varphi(u)$ be an nonnegative even function
such that it is not identically zero on the real line, if it is indefinitely
differentiable and for each $n\in\mathbb{Z}^{+}$, 
\begin{equation}
\lim_{u\to\infty}\frac{\log\left|\varphi^{(n)}(u)\right|}{u}=-\infty.\label{eq:2.1}
\end{equation}
Additionally, the even entire function 
\begin{equation}
\Phi(z)=\int_{-\infty}^{\infty}\varphi(u)e^{iuz}du=2\int_{0}^{\infty}\varphi(u)\cos(zu)du\label{eq:2.2}
\end{equation}
is of genus at most $1$. Then the inequality 
\begin{equation}
\frac{\partial^{2}}{\partial y^{2}}\left|\Phi\left(x+iy\right)\right|^{2}\ge0,\quad x,\,y\in\mathbb{R},\label{eq:2.3}
\end{equation}
is necessary and sufficient for $\Phi(z)$ to have only real zeros. \end{thm}
\begin{lem}
\label{lem:2} Let $f(z),\ h(z)$ be two entire functions of genus
at most $1$, then so is $f(z)h(z)$.\end{lem}
\begin{proof}
Let $\rho_{f}$ and $\rho_{h}$ be the orders of $f$ and $h$ respectively.
Since they are of genus at most $1$, $0\le\left\lfloor \rho_{f}\right\rfloor ,\,\left\lfloor \rho_{h}\right\rfloor \le1$,
then we must have $0\le\rho_{f},\rho_{h}<2.$ By (\ref{eq:1.2}) we
have 
\begin{equation}
\frac{\log\left(\log\left\Vert f\right\Vert _{\infty,r}\right)}{\log(r)}\le a,\quad\frac{\log\left(\log\left\Vert h\right\Vert _{\infty,r}\right)}{\log(r)}\le b,\label{eq:2.4}
\end{equation}
as $r\to\infty$ where $0<a,\ b<2$. Thus,
\begin{equation}
\log\left\Vert f\right\Vert _{\infty,r}\le r^{a},\quad\log\left\Vert h\right\Vert _{\infty,r}\le r^{b}\label{eq:2.5}
\end{equation}
as $r\to\infty$. From 
\begin{equation}
\log\left(\left\Vert fh\right\Vert _{\infty,r}\right)\le\log\left(\left\Vert f\right\Vert _{\infty,r}\cdot\left\Vert h\right\Vert _{\infty,r}\right)\le\log\left(\left\Vert f\right\Vert _{\infty,r}\right)+\log\left(\left\Vert h\right\Vert _{\infty,r}\right)\label{eq:2.6}
\end{equation}
we get
\begin{equation}
\log\left(\left\Vert fh\right\Vert _{\infty,r}\right)\le r^{a}+r^{b}\le2r^{\max\left\{ a,b\right\} }\label{eq:2.7}
\end{equation}
for $r\to\infty$. Then,
\begin{equation}
\log\left(\log\left(\left\Vert fh\right\Vert _{\infty,r}\right)\right)\le\log2+\max\left\{ a,b\right\} \log r,\label{eq:2.8}
\end{equation}
and
\begin{equation}
\limsup_{r\to\infty}\frac{\log\left(\log\left(\left\Vert fh\right\Vert _{\infty,r}\right)\right)}{\log r}\le\max\left\{ a,b\right\} <2.\label{eq:2.9}
\end{equation}
Hence the order of $fh$ is strictly less $2$, thus the genus of
$fh$ is less than $\left\lfloor \max\left\{ a,b\right\} \right\rfloor \le1$.\end{proof}
\begin{lem}
\label{lem:3}Let $\varphi(u)$ be an even nonnegative measurable
function that is not identically zero. If there exists some $d,\,\delta>0$
such that $\varphi(u)=\mathcal{O}\left(\exp\left(-de^{\delta\left|u\right|}\right)\right)$
as $u\to\pm\infty$. Then the entire function defined in (\ref{eq:2.2})
is of genus at most $1$\end{lem}
\begin{proof}
By the assumptions we have
\begin{equation}
b_{2n}=\int_{-\infty}^{\infty}\varphi(u)u^{2n}du>0,\quad n\in\mathbb{Z}^{+}\label{eq:2.10}
\end{equation}
and
\begin{equation}
\begin{aligned} & \Phi(z)=2\int_{0}^{\infty}\varphi(u)\cos(zu)du\\
 & =2\int_{0}^{\infty}\varphi(u)\left(\sum_{n=0}^{\infty}\frac{\left(-z^{2}u^{2}\right)^{n}}{(2n)!}\right)du=\sum_{n=0}^{\infty}\frac{\left(-z^{2}\right)^{n}}{(2n)!}b_{2n}.
\end{aligned}
\label{eq:2.11}
\end{equation}
Then,
\begin{equation}
\sup_{|z|\le r}\left|\Phi(z)\right|=\sup_{|z|\le r}\left|\sum_{n=0}^{\infty}\frac{\left(-z^{2}\right)^{n}}{(2n)!}b_{2n}\right|\le\sum_{n=0}^{\infty}\frac{r^{2n}}{(2n)!}b_{2n}=\Phi(ir).\label{eq:2.12}
\end{equation}
By (\ref{eq:2.1}) we have
\begin{equation}
\begin{aligned} & \Phi(ir)=2\int_{0}^{\infty}\varphi(u)\cosh(ru)du=\int_{0}^{\infty}\mathcal{O}\left(\exp\left(-de^{\delta u}+ru\right)\right)du\\
 & =\mathcal{O}\left(\int_{0}^{\infty}\exp\left(-de^{\delta u}+ru\right)du\right)=\mathcal{O}\left(\frac{1}{d^{\frac{r}{\delta}}}\int_{d}^{\infty}\exp\left(-y\right)y^{\frac{r}{\delta}-1}du\right)\\
 & =\mathcal{O}\left(\frac{1}{d^{\frac{r}{\delta}}}\int_{0}^{\infty}\exp\left(-y\right)y^{\frac{r}{\delta}-1}du\right)=\mathcal{O}\left(\frac{\Gamma\left(\frac{r}{\delta}\right)}{d^{\frac{r}{\delta}}}\right).
\end{aligned}
\label{eq:2.13}
\end{equation}
By Stirling's formula we have \cite{Andrews}
\begin{equation}
\begin{aligned} & \log\left(\Phi(ir)\right)=-\frac{r}{\delta}\log d+\log\Gamma\left(\frac{r}{\delta}\right)+\mathcal{O}(1)\\
 & =\mathcal{O}\left(r\log r\right)
\end{aligned}
\label{eq:2.14}
\end{equation}
as $r\to+\infty$. Hence,
\begin{equation}
\begin{aligned}\rho= & \limsup_{r\to\infty}\frac{\log\left(\log\Phi(ir)\right)}{\log(r)}\le1,\end{aligned}
\label{eq:2.15}
\end{equation}
which shows that the order $\rho$ of $\Phi(z)$ is no higher than
$1$, thus $\Phi(z)$ is of genus $0$ or $1$. 
\end{proof}
Let 
\begin{equation}
t(u)=\begin{cases}
\left(1-u^{2}\right),\  & -1<u<1,\\
0,\  & \left|u\right|\ge1.
\end{cases}\label{eq:2.16}
\end{equation}
Clearly, 
\begin{equation}
\frac{d^{2}t(u)}{du^{2}}=\begin{cases}
-2,\  & -1<u<1\\
0,\  & \left|u\right|>1
\end{cases}\label{eq:2.17}
\end{equation}
 and \cite{Andrews,Erdelyi,Gasper1,Ismail,Watson}
\begin{equation}
\sqrt{8\pi}z^{-\frac{3}{2}}J_{\frac{3}{2}}(z)=\int_{-1}^{1}e^{izu}\left(1-u^{2}\right)du.\label{eq:2.18}
\end{equation}

\begin{lem}
\label{lem:4} Let $\varphi(u)$ be nonnegative, even, and indefinitely
differentiable function that is not identically zero. Additionally,
we assume that for each $n\in\mathbb{Z}^{+}$ there exists some $d>0$
such that 
\begin{equation}
\varphi^{(n)}(u)=\mathcal{O}\left(\exp\left(-de^{\delta\left|u\right|}\right)\right)\label{eq:2.19}
\end{equation}
as $u\to\pm\infty$. Let us define 
\begin{equation}
\psi(u)=\int_{\mathbb{R}}\varphi(u-y)t(y)dy,\quad u\in\mathbb{R}.\label{eq:2.20}
\end{equation}
Then, $\psi(u)$ is nonnegative, even, indefinitely differentiable
function that is not identically zero and 
\begin{equation}
\psi^{(2)}(u)\le0,\quad u\in\mathbb{R}.\label{eq:2.21}
\end{equation}
Furthermore, for each $n\in\mathbb{Z}^{+}$ we have 
\begin{equation}
\psi^{(n)}(u)=\mathcal{O}\left(\exp\left(-d_{\delta}e^{\delta\left|u\right|}\right)\right)\label{eq:2.22}
\end{equation}
as $u\to\pm\infty$, where $d_{\delta}=\frac{d}{e^{\delta}}$. \end{lem}
\begin{proof}
Since for $u\in\mathbb{R}$ and $|y|\le1$ we have
\begin{equation}
|u-y|\ge|u|-1,\quad e^{\delta\left|u-y\right|}\ge e^{-\delta}e^{\delta\left|u\right|},\quad\exp\left(-de^{\delta\left|u-y\right|}\right)\le\exp\left(-\frac{d}{e^{\delta}}e^{\delta\left|u\right|}\right).\label{eq:2.23}
\end{equation}
Then by the asymptotic behavior of $\varphi^{(n)}(u)$, we have
\begin{equation}
\int_{\mathbb{R}}\varphi^{(n)}(u-y)t(y)dy=\mathcal{O}\left(\exp\left(-\frac{d}{e^{\delta}}e^{\delta\left|u\right|}\right)\int_{-1}^{1}t(y)dy\right)=\mathcal{O}\left(\exp\left(-d_{\delta}e^{\delta\left|u\right|}\right)\right),\label{eq:2.24}
\end{equation}
where $d_{\delta}=\frac{d}{e^{\delta}}$. Thus the integral $\int_{\mathbb{R}}\varphi^{(n)}(u-y)t(y)dy$
converges absolutely and uniformly for each $n\in\mathbb{Z}^{+}$,
hence 
\begin{equation}
\psi^{(n)}(u)=\int_{\mathbb{R}}\varphi^{(n)}(u-y)t(y)dy=\mathcal{O}\left(\exp\left(-d_{\delta}e^{\delta\left|u\right|}\right)\right)\label{eq:2.25}
\end{equation}
as $u\to\pm\infty$. Clearly, $\psi(u)$ is nonnegative and not identically
zero if $\varphi(u)$ is such, and for each $u\in\mathbb{R}$ we have
\begin{equation}
\psi(u)=\int_{\mathbb{R}}\varphi(u+y)t(-y)dy=\int_{\mathbb{R}}\varphi(-u-y)t(y)dy=\psi(-u).\label{eq:2.26}
\end{equation}
 Since
\begin{equation}
\psi(u)=\int_{\mathbb{R}}\varphi(u-y)t(y)dy=\int_{\mathbb{R}}\varphi(y)t(u-y)dy,\label{eq:2.27}
\end{equation}
 then,
\begin{equation}
\begin{aligned}\psi^{(2)}\left(u\right) & =\int_{\mathbb{R}}\varphi(y)t^{(2)}\left(x-y\right)dy=-2\int_{|x-y|\le1}\varphi(y)dy\le0.\end{aligned}
\label{eq:2.28}
\end{equation}
\end{proof}
\begin{lem}
\label{lem:5} Let 
\begin{equation}
F(x)=\int_{0}^{\infty}G(y)\cos xydy,\label{eq:2.29}
\end{equation}
where $G(y)\in\mathcal{C}^{2}\left([0,\infty\right)$ and 
\begin{equation}
\int_{0}^{\infty}\left\{ \left|G(y)\right|+\left|G'(y)\right|+\left|G''(y)\right|\right\} dy<\infty.\label{eq:2.30}
\end{equation}
If $G'(0)=0$ and $G''(y)\le0$ for $y\ge0$ almost everywhere, then
$F(x)\ge0$ for all $x\ge0$. In particular, if $G(y)$ is even and
$G''(y)\le0$ for $y\ge0$ almost everywhere, then $F(x)\ge0$ for
all $x\ge0$. \end{lem}
\begin{proof}
Observe that for $G''(y)\le0$ and $G'(0)=0$ we have
\begin{equation}
\begin{aligned} & x^{2}F(x)=x\int_{0}^{\infty}G(y)d\sin xy=-x\int_{0}^{\infty}G'(y)\sin xydy\\
 & =\int_{0}^{\infty}G'(y)d(1+\cos xy)=-2\int_{0}^{\infty}G''(y)\cos^{2}\left(\frac{xy}{2}\right)dy\ge0.
\end{aligned}
\label{eq:2.31}
\end{equation}
 In the case that $G(y)$ is even, then $G'(y)$ is odd and $G'(0)=0$.
\end{proof}

Now we prove our main result Theorem \ref{thm:gasper}:
\begin{proof}
It is known that the even entire function $\sqrt{8\pi}z^{-\frac{3}{2}}J_{\frac{3}{2}}(z)$
is of genus $1$ and it has only real zeros, \cite{Andrews,Gasper1,Gasper2,Ismail,Watson}. 

Let 
\begin{equation}
\Psi(z)=\sqrt{8\pi}z^{-\frac{3}{2}}J_{\frac{3}{2}}(z)\Phi(z)=\int_{\mathbb{R}}\psi(u)e^{izu}\,du,\label{eq:2.32}
\end{equation}
by (\ref{eq:2.18}) we have
\begin{equation}
\Psi(z)=\int_{\mathbb{R}}\psi(u)e^{izu}\,du,\label{eq:2.33}
\end{equation}
where $\psi(u)$ is defined by (\ref{eq:2.20}). By Lemmas \ref{lem:2},
\ref{lem:3} and \ref{lem:4} we know that in addition to $\Psi(z)$
and $\psi(u)$ satisfy all the requirements of Theorem \ref{thm:jensen}
except (\ref{eq:2.3}). Furthermore, we also have $\psi^{(2)}(u)\le0,\ u\in\mathbb{R}$. 

From

\begin{equation}
\overline{\Psi(z)}=\int_{-\infty}^{\infty}\psi(u)e^{-i\overline{z}u}du=\int_{-\infty}^{\infty}\psi(u)e^{i\overline{z}u}du=\Psi(\overline{z}),\label{eq:2.34}
\end{equation}
we get 
\begin{equation}
\begin{aligned} & 2\left|\Psi(z)\right|^{2}=2\Psi(z)\Psi(\overline{z})=2\int_{\mathbb{R}^{2}}\psi(u)\psi(v)e^{-y(u-v)}e^{ix(u+v)}dudv\\
 & =\int_{\mathbb{R}^{2}}\psi\left(\frac{\alpha+\beta}{2}\right)\psi\left(\frac{\beta-\alpha}{2}\right)e^{-\alpha y}e^{i\beta x}d\alpha d\beta\\
 & =\int_{\mathbb{R}^{2}}\psi\left(\frac{\alpha+\beta}{2}\right)\psi\left(\frac{\alpha-\beta}{2}\right)e^{-\alpha y}e^{i\beta x}d\alpha d\beta\\
 & =\int_{\mathbb{R}^{2}}\psi\left(\frac{-\alpha+\beta}{2}\right)\psi\left(\frac{-\alpha-\beta}{2}\right)e^{\alpha y}e^{i\beta x}d\alpha d\beta\\
 & =\int_{\mathbb{R}^{2}}\psi\left(\frac{\alpha+\beta}{2}\right)\psi\left(\frac{\alpha-\beta}{2}\right)e^{\alpha y}e^{i\beta x}d\alpha d\beta\\
 & =\int_{\mathbb{R}^{2}}\psi\left(\frac{\alpha+\beta}{2}\right)\psi\left(\frac{\alpha-\beta}{2}\right)\cosh(y\alpha)e^{i\beta x}d\alpha d\beta\\
 & =\int_{\mathbb{R}^{2}}\psi\left(\frac{\alpha+\beta}{2}\right)\psi\left(\frac{\alpha-\beta}{2}\right)\cosh(y\alpha)e^{-i\beta x}d\alpha d\beta\\
 & =2\int_{0}^{\infty}\left\{ \int_{\mathbb{R}}\psi\left(\frac{\alpha+\beta}{2}\right)\psi\left(\frac{\alpha-\beta}{2}\right)\cosh(y\alpha)d\alpha\right\} \cos(\beta x)d\alpha d\beta.
\end{aligned}
\label{eq:2.35}
\end{equation}
Then
\begin{equation}
\begin{aligned} & \frac{\partial^{2}}{\partial y^{2}}\left|\Psi(z)\right|^{2}=\int_{0}^{\infty}\left\{ \int_{\mathbb{R}}\psi\left(\frac{\alpha+\beta}{2}\right)\psi\left(\frac{\alpha-\beta}{2}\right)\alpha^{2}\cosh(y\alpha)d\alpha\right\} \cos(\beta x)d\beta.\end{aligned}
\label{eq:2.36}
\end{equation}
Let 
\begin{equation}
G(\beta;y)=\int_{\mathbb{R}}\psi\left(\frac{\alpha+\beta}{2}\right)\psi\left(\frac{\alpha-\beta}{2}\right)\alpha^{2}\cosh(y\alpha)d\alpha.\label{eq:2.37}
\end{equation}
From
\begin{equation}
\begin{aligned} & \partial_{\alpha}^{2}\left\{ \psi\left(\frac{\alpha+\beta}{2}\right)\psi\left(\frac{\alpha-\beta}{2}\right)\right\} +\partial_{\beta}^{2}\left\{ \psi\left(\frac{\alpha+\beta}{2}\right)\psi\left(\frac{\alpha-\beta}{2}\right)\right\} \\
 & =\frac{1}{2}\left\{ \psi^{(2)}\left(\frac{\alpha+\beta}{2}\right)\psi\left(\frac{\alpha-\beta}{2}\right)+\psi\left(\frac{\alpha+\beta}{2}\right)\psi^{(2)}\left(\frac{\alpha-\beta}{2}\right)\right\} 
\end{aligned}
\label{eq:2.38}
\end{equation}
we get
\begin{equation}
\begin{aligned} & \partial_{\beta}^{2}G(\beta;y)=\partial_{\beta}^{2}\int_{\mathbb{R}}\psi\left(\frac{\alpha+\beta}{2}\right)\psi\left(\frac{\alpha-\beta}{2}\right)\alpha^{2}\cosh(y\alpha)d\alpha\\
 & =\int_{\mathbb{R}}\partial_{\beta}^{2}\left\{ \psi\left(\frac{\alpha+\beta}{2}\right)\psi\left(\frac{\alpha-\beta}{2}\right)\right\} \alpha^{2}\cosh(y\alpha)d\alpha\\
 & =\frac{1}{2}\int_{\mathbb{R}}\left\{ \psi^{(2)}\left(\frac{\alpha+\beta}{2}\right)\psi\left(\frac{\alpha-\beta}{2}\right)+\psi\left(\frac{\alpha+\beta}{2}\right)\psi^{(2)}\left(\frac{\alpha-\beta}{2}\right)\right\} \alpha^{2}\cosh(y\alpha)d\alpha\\
 & -\int_{\mathbb{R}}\partial_{\alpha}^{2}\left\{ \psi\left(\frac{\alpha+\beta}{2}\right)\psi\left(\frac{\alpha-\beta}{2}\right)\right\} \alpha^{2}\cosh(y\alpha)d\alpha\\
 & =\frac{1}{2}\int_{\mathbb{R}}\left\{ \psi^{(2)}\left(\frac{\alpha+\beta}{2}\right)\psi\left(\frac{\alpha-\beta}{2}\right)+\psi\left(\frac{\alpha+\beta}{2}\right)\psi^{(2)}\left(\frac{\alpha-\beta}{2}\right)\right\} \alpha^{2}\cosh(y\alpha)d\alpha\\
 & -\int_{\mathbb{R}}\psi\left(\frac{\alpha+\beta}{2}\right)\psi\left(\frac{\alpha-\beta}{2}\right)\left((2+a^{2}y^{2})\cosh(ay)+4(ay)\sinh(ay)\right)d\alpha\le0.
\end{aligned}
\label{eq:2.39}
\end{equation}
Then by Lemma \ref{lem:5} we have proved that $\frac{\partial^{2}}{\partial y^{2}}\left|\Psi(z)\right|^{2}\ge0$.
Then by Jensen's result Theorem \ref{thm:jensen} we have shown that
$\Psi(z)=\sqrt{8\pi}z^{-\frac{3}{2}}J_{\frac{3}{2}}(z)\Phi(z)$ has
only real zeros. Since $\sqrt{8\pi}z^{-\frac{3}{2}}J_{\frac{3}{2}}(z)$
has only real zeros, therefore, $\Phi(z)$ has only real zeros.\end{proof}


\begin{thebibliography}{10}
\bibitem{Gasper2} G. Gasper, Using integrals of squares to prove
that certain real-valued special functions to prove that the P\'{o}lya
$\Xi^{*}\left(z\right)$ function, the functions $K_{iz}\left(a\right),\ a>0$
and some other entire functions have only real zeros, dedicated to
Dan Waterman on the occasion of his 80th birthday.
\end{thebibliography}
\end{document}